\documentclass[11pt]{article}

\usepackage{amssymb}
\usepackage{amsmath}
\usepackage{theorem}
\usepackage{epsfig}
\usepackage{verbatim}
\usepackage{graphicx}

\newtheorem{thm}{Theorem}[section]

\newtheorem{prop}[thm]{Proposition}

\newtheorem{lemma}[thm]{Lemma}

\newcommand{\R}{\Bbb{R}}

\newcommand{\N}{\Bbb{N}}
\newcommand{\T}{\mathbb{T}}
\newcommand{\D}{\displaystyle}

\newcommand{\di}{{\rm div}\thinspace}
\newcommand{\curl}{{\rm curl}\thinspace}

\newcommand{\grad}{\nabla}

\newcommand{\dx}{\partial_x}

\newcommand{\dxu}{\partial_{x_1}}
\newcommand{\dxd}{\partial_{x_2}}
\newcommand{\dxt}{\partial_{x_3}}

\newcommand{\la}{\Lambda}

\newcommand{\al}{\alpha}

\begin{document}

\author{Diego C\'ordoba and Francisco Gancedo}
\title{A maximum principle for the Muskat problem\\ for fluids with different densities}
\date{ }

\maketitle

\begin{abstract}
We consider the fluid interface problem given by two incompressible fluids with different densities evolving by Darcy's law.
This scenario is known as the Muskat problem for fluids with the same viscosities, being in two dimensions mathematically analogous to the two-phase Hele-Shaw cell. We prove in the stable case (the denser fluid is below) a maximum principle for the $L^\infty$ norm of the free boundary.
\end{abstract}

\maketitle

%%%%%%%%%%%%%%%%%%%%%%%%%%%%%%%%%%%%%%%%%%%%%%%%%%%%%%%%%%%%%%%%%%%%%%%%%%%5
\section{Introduction}

The Muskat problem models the fluid interface problem given by two fluids in a
porous medium with different characteristics. The problem was proposed by Muskat (see \cite{Muskat})
in a study about the encroachment of water into oil in a porous medium.
In this phenomena, Darcy's law is used to govern the dynamics of the different fluids \cite{bear}. This law is given by 
the following formula:
\begin{equation*}\label{dl}
\frac{\mu}{\kappa} v=-\nabla p-(0,0,\mathrm{g}\,\rho),
\end{equation*}
where $v$ is the velocity of the fluid, $p$ is the pressure, $\mu$ is the dynamic viscosity,
$\kappa$ is the permeability of the isotropic medium, $\rho$ is the liquid density and $\mathrm{g}$ is the
acceleration due to gravity. The fluid interface is given between incompressible flows with different viscosities, $\mu_1$, $\mu_2$, 
and densities, $\rho_1$, $\rho_2$.

Saffman and Taylor \cite{S-T} considered this problem in a study of the dynamics of the interface between two fluids with different
viscosities and densities in a Hele--Shaw cell. In this physical scenario (see \cite{H-S})
the fluid is trapped between two fixed parallel plates, that are close enough together, so
that the fluid essentially only moves in two directions. The mean velocity is described by
\begin{equation*}
\frac{12\mu}{b^2} v=-\nabla p-(0,\mathrm{g}\,\rho),
\end{equation*}
where $b$ is the distance between the plates. Darcy's law, in two dimensions, and the above formula become analogous if we consider the permeability of the medium 
$\kappa$ equal to the constant $b^2/12$.

The Muskat problem and the two--phase Hele--Shaw flow have been extensively considered (see \cite{Peter} and \cite{Hou} and the references therein). These free boundary problems can be modeled with surface tension \cite{ES} using the Laplace--Young condition. In this case there is a jump of discontinuity in the pressure of the fluids across the interface proportional to the local curvature of the free boundary. In this work we will consider the case without surface tension, so that the pressures are equal on the interface. 

This scenario is considered by Siegel, Caflisch and Howison in \cite{SCH}, where they show ill-posedness
in an unstable case and global-in-time existence of small initial data in a stable case. They describe the two-dimensional dynamics of the incompressible flow as follows:
$$
v=-a\nabla p-(0,V),
$$
where $a$ takes two positive constant values
$$
a_i=\frac{b^2}{12\mu_i},\qquad\mbox{for}\qquad i=1,2,
$$ on each fluid, and $V$ is a constant. With our notation, this case is equivalent to consider in the two-dimensional problem $$\frac{\mu_i}{\kappa}=\frac{1}{a_i},\qquad\mbox{and}\qquad \mathrm{g}\,\rho_i=\frac{V}{a_i},$$ for $i=1,2$. The results rely on the assumption that the 
Atwood number 
$$
A_\mu=\frac{\mu_1-\mu_2}{\mu_1+\mu_2},
$$ is nonzero, and therefore there is a jump of viscosities. 

In the same year, Ambrose \cite{Ambrose} studies the 2-D problem with initial data fulfilling
$$
(\rho_2-\rho_1)\mathrm{g}\,cos(\theta(\alpha,0))+2A_\mu U(\alpha,0)>0,
$$
where $\theta$ is the angle that the tangent to the curve forms with the horizontal and $U$
is the normal velocity (given by the Birkhoff-Rott integral). In this work he considers the arclength and the 
tangent angle formulation used by Hou, Lowengrub and Shelley \cite{Hou} to get energy 
estimates for the free boundary assuming the geometric condition 
\begin{equation}\label{ambrose}
\frac{(x(\alpha,t)-x(\alpha',t))^2+(y(\alpha,t)-y(\alpha',t))^2}{(\al-\al')^2}>0,
\end{equation}
locally in time, where the curve $(x(\al,t),y(\al,t))$ is the interface. We show in \cite{Y} that this is not enough to get local-existence in this kind of contour dynamics equations, since a regular interface could touch itself with order infinity and without satisfying the inequality \eqref{ambrose}.

We will study the fluid interface given by a jump of densities, so that $A_\mu=0$. Therefore we can take $\kappa=b^2/12=1$, and $\mu_1=\mu_2=\mathrm{g}=1$ without loss of generality. This case can consider, among others, the dynamics of moist and dry regions in porous media. This scenario was considered 
by Dombre, Pumir and Siggia \cite{DPS}, but in a different context. They treated the interface dynamics for
convection in porous media where, using our notation, the density plays the roll of the temperature. They studied the unstable case, namely when the denser fluid (or the fluid with larger temperature) is above. They present meromorphic initial conditions with complex poles, and study the dynamics of these critical points.

It is well-known \cite{Hou} that for these contour dynamics systems, the velocity in the tangential direction does not alter the shape of the interface. If we change the tangential component of the velocity, we only change the parametrization. In \cite{DY}
we used this property to parameterize the interface as a function $(x,f(x,t))$, getting the following equations:
\begin{align}
\begin{split}\label{ecintro}
\D f_t (x,t)&=\frac{\rho_2-\rho_1}{4\pi}PV\int_{\R^2}\frac{(\grad f(x,t)-\grad
f(x-y,t))\cdot y}{[|y|^2+(f(x,t)-f(x-y,t))^2]^{3/2}}dy,\\
f(x,0)&=f_0(x),\quad x\in\R^2,
\end{split}
\end{align}
for a two-dimensional interface, and
\begin{align}
\begin{split}\label{ec1dintro}
\D f_t(\al,t)&=\frac{\rho_2-\rho_1}{2\pi}PV\int_{\R}\frac{(\dx f(\al,t)-\dx
f(\al-\beta,t))\beta}{\beta^2+(f(\al,t)-f(\al-\beta,t))^2}d\beta,\\
f(\al,0)&=f_0(\al),\quad \al\in\R,
\end{split}
\end{align}
for a one-dimensional interface. We point out that with these formulations the condition \eqref{ambrose} is satisfied locally in time if local-existence for the systems is reached. This avoids a kind of singularity in the fluid when the interface collapses (see \cite{CFMR} for example). Also we prove that when the denser fluid is below the other fluid, $\rho_2>\rho_1$, the problem is well-posed given local-existence and uniqueness for the systems \eqref{ecintro} and \eqref{ec1dintro}. When the less dense fluid is below, $\rho_2<\rho_1$, we prove ill-posedness showing that the equations \eqref{ecintro} and \eqref{ec1dintro} are ill-posed. We get this result using global solutions of \eqref{ec1dintro} in the stable case, $\rho_2>\rho_1$, for small initial data in a similar way as in \cite{SCH}. 

If we neglect the terms of order two in \eqref{ecintro}, the linearized equation is obtained. It reads
\begin{align}
\begin{split}\label{ecl}
&f_t=\frac{\rho_1-\rho_2}{2}(R_1\dxu f+R_2\dxd f)=\frac{\rho_1-\rho_2}{2}\la f,\\
&f(x,0)=f_0(x),
\end{split}
\end{align}
where $R_1$ and $R_2$ are the Riesz transforms (see \cite{St3}) and the operator $\la f$ is given by the Fourier transform $\widehat{\la f}(\xi)=|\xi|\widehat{f}(\xi)$. In the stable case $\rho_1<\rho_2$ (the greater density is below), the linear equation is dissipative and is clear that the following maximum principle is reached
$$
\|f\|_{L^{\infty}}(t)\leq \|f_0\|_{L^{\infty}}.
$$
In this work we get, in section 3 and 4, this estimate for the nonlinear systems \eqref{ecintro} and \eqref{ec1dintro}. To this end, we follow the evolution of the maximum of the absolute value of $f(x,t)$. This technique was used for one of the authors in \cite{cor2} in a family of dissipative transport equations for incompressible fluids. Also we would like to cite the work of A. Constantin and J. Echer where they study the shallow water equation in the same way. By a similar approach, in section 5 we obtain a global bound on the derivative for small initial data.

%The free boundary problems given by fluids with different densities are been intensely studied.
%Notice the classical paper of Taylor \cite{Taylor} and the works of Wu \cite{Wu} and \cite{Wu2}
%where the full water wave problem is solved considering the water with positive density and the air
%with zero density. A study of the two-dimensional case can be found in \cite{AM} due to Ambrose and
%Masmoudi.

%%%%%%%%%%%%%%%%%%%%%%%%%%%%%%%%%%%%%%%%%%%%%%%%%%%%%%%%%%%%%%%%%%%%%%%%%%%%%%%%%%%%%%%%%%%%%%%%

%%%%%%%%%%%%%%%%%%%%%%%%%%%%%%%%%%%%%%%%%%%%%%%%%%%%%%%%%%%%%%%%%%%%%%%%%%%%%%%%%%%%%%%%%%%%%%%%

%%%%%%%%%%%%%%%%%%%%%%%%%%%%%%%%%%%%%%%%%%%%%%%%%%%%%%%%%%%%%%%%%%%%%%%%%%%%%%%%%%%%%%%%%%%%%%%%

\section{Parameterizing the interface in terms of a function}

%%%%%%%%%%%%%%%%%%%%%%%%%%%%%%%%%%%%%%%%%%%%%%%%%%%%%%%%%%%%%%%%%%%%%%%%%%%%%%%%%%%%%%%%%%%%%%%%

%%%%%%%%%%%%%%%%%%%%%%%%%%%%%%%%%%%%%%%%%%%%%%%%%%%%%%%%%%%%%%%%%%%%%%%%%%%%%%%%%%%%%%%%%%%%%%%%

%%%%%%%%%%%%%%%%%%%%%%%%%%%%%%%%%%%%%%%%%%%%%%%%%%%%%%%%%%%%%%%%%%%%%%%%%%%%%%%%%%%%%%%%%%%%%%%%

In this section we briefly explain how to parameterize the free boundary in terms of a function (see \cite{DY} for more details). This way of writing this nonlocal equation is crucial to check the evolution of the maximum of the absolute value of the function to obtain the maximum principle. 

In our case, Darcy's law can be written as follow:
\begin{equation}\label{dl2}
v(x_1,x_2,x_3,t)=-\nabla p(x_1,x_2,x_3,t)-(0,0,\rho(x_1,x_2,x_3,t)),
\end{equation}
where $(x_1,x_2,x_3)\in\R^3$ are the spatial variables and $t\geq 0$ denotes the time. Here $\rho$ is defined by
$$
\rho(x_1,x_2,x_3,t)=\left\{\begin{array}{cl}
                    \rho_1\quad\mbox{in}&\Omega_1(t)\\
                    \rho_2\quad\mbox{in}&\Omega_2(t),
                 \end{array}\right. $$ 
with $\rho_1,\rho_2\geq 0$ constants and $\rho_1\neq\rho_2$. The set $\Omega_1(t)$ is given by $$\Omega_1(t)=\{x_3> f(x_1,x_2,t)\},$$ and $$\Omega_2(t)=\{x_3< f(x_1,x_2,t)\},$$
being $f(x_1,x_2,t)$ the fluid interface. If we apply the $\curl$ operator to Darcy's law twice the pressure disappears. Considering the incompressibility of the fluid, we have $\curl\curl v=-\Delta v$, and we can give the velocity in terms of the density as follows: 
\begin{equation}\label{sior}
v=(\dxu\Delta^{-1}\dxt\rho,\dxd\Delta^{-1}\dxt\rho,-\dxu\Delta^{-1}\dxu\rho-\dxd\Delta^{-1}\dxd\rho).
\end{equation}
The density $\rho$ has a jump of discontinuity on the free boundary, therefore the gradient of it function is given by a Dirac distribution $\delta$ as follows:
\begin{equation}\label{gradprho}
\grad\rho=(\rho_2-\rho_1)(\dxu f(x_1,x_2,t),\dxd f(x_1,x_2,t),-1)\delta(x_3-f(x_1,x_2,t)),
\end{equation} Using the kernels for $\dxu\Delta^{-1}$ and $\dxd\Delta^{-1}$ we obtain
\begin{align}\label{pesao}
\begin{split}
v(x_1,x_2,x_3,t)=&-\D\frac{\rho_2-\rho_1}{4\pi}PV\int_{\R^2}\frac{(y_1,y_2,\grad f(x-y,t)\cdot y)}
{[|y|^2+(x_3-f(x-y,t))^2]^{3/2}}dy,
\end{split}
\end{align} where we note $x=(x_1,x_2)$, $y=(y_1,y_2)$. In \eqref{pesao} $x_3\neq f(x,t)$ and the principal value is taken at infinity (see \cite{St3}). The vorticity is at the same level than the gradient of the density, so it is given by a delta function. This forces the velocity to have a discontinuity on the free boundary. Just checking the incompressibility of the fluid in the sense of the distributions, we obtain that this discontinuity is in the tangential direction, so that it does not affect the shape of the interface (see \cite{DY}). Ignoring the tangential terms we obtain that the velocity on the free boundary is given by
\begin{align}\label{fvpv}
\begin{split}
v(x,f(x,t),t)=&-\D\frac{\rho_2-\rho_1}{4\pi}PV\int_{\R^2}\frac{(y_1,y_2,\grad f(x-y,t)\cdot y)}
{[|y|^2+(f(x,t)-f(x-y,t))^2]^{3/2}}dy.
\end{split}
\end{align} If we want to parameterice the free boundary in terms of a function, we need to get the velocity
$v=(v_1,v_2,v_3)$ satisfying $v_1=v_2=0$, since otherwise the points on the plane are not fixed and they depend on time. If we add the 
following tangential terms to \eqref{fvpv}:
$$
\frac{\rho_2-\rho_1}{4\pi}PV\int_{\R^2}\frac{y_1}{[|y|^2+(f(x,t)-f(x-y,t))^2]^{3/2}}dy(1,0,\dxu
f(x,t)),
$$
$$
\frac{\rho_2-\rho_1}{4\pi}PV\int_{\R^2}\frac{y_2}{[|y|^2+(f(x,t)-f(x-y,t))^2]^{3/2}}dy(0,1,\dxd
f(x,t)),
$$
we do not alter the interface and this follows:
\begin{align}\label{vn}
v(x,f(x,t),t)&=\frac{\rho_2-\rho_1}{4\pi}(0,0,PV\int_{\R^2}\frac{(\grad f(x,t)-\grad f(x-y,t))\cdot
y}{[|y|^2+(f(x,t)-f(x-y,t))^2]^{3/2}}dy).
\end{align}
Finally the contour equation given by
\begin{align}
\begin{split}\label{ec}
\D f_t (x,t)&=\frac{\rho_2-\rho_1}{4\pi}PV\int_{\R^2}\frac{(\grad f(x,t)-\grad
f(x-y,t))\cdot y}{[|y|^2+(f(x,t)-f(x-y,t))^2]^{3/2}}dy,\\
f(x,0)&=f_0(x).
\end{split}
\end{align}
This formula works for periodic interface and for a free boundary near planar at infinity. In both cases it presents an principal value only at infinity. If we suppose that the function $f(x,t)$ only depends on $x_1$, integrating in $x_2$, the contour equation in the 2-D case is reached
\begin{align}
\begin{split}\label{ec1d}
\D f_t(x,t)&=\frac{\rho_2-\rho_1}{2\pi}PV\int_{\R}\frac{(\dx f(x,t)-\dx
f(x-\al,t))\al}{\al^2+(f(x,t)-f(x-\al,t))^2}d\al,\\
f(x,0)&=f_0(x);\quad x\in\R.
\end{split}
\end{align}
This equation can be obtained in a similar way that \eqref{ec} using the stream function. 
We check in \cite{DY} that as long as this equation is satisfied we obtain weak solutions of the system
\begin{equation}\label{fsystem}
\begin{array}{l}
\rho_t+v\cdot\grad \rho=0,\\
\\
v=-\grad p-(0,0,\rho),\quad \di v=0. 
\end{array}
\end{equation}

%%%%%%%%%%%%%%%%%%%%%%%%%%%%%%%%%%%%%%%%%%%%%%%%%%%%%%%%%%%%%%%%%%%%%%%%%%%%%%%%%%%%%%%%%%%%%%%%

%%%%%%%%%%%%%%%%%%%%%%%%%%%%%%%%%%%%%%%%%%%%%%%%%%%%%%%%%%%%%%%%%%%%%%%%%%%%%%%%%%%%%%%%%%%%%%%%

%%%%%%%%%%%%%%%%%%%%%%%%%%%%%%%%%%%%%%%%%%%%%%%%%%%%%%%%%%%%%%%%%%%%%%%%%%%%%%%%%%%%%%%%%%%%%%%%

\section{Two dimensional case (1-D interface)}

%%%%%%%%%%%%%%%%%%%%%%%%%%%%%%%%%%%%%%%%%%%%%%%%%%%%%%%%%%%%%%%%%%%%%%%%%%%%%%%%%%%%%%%%%%%%%%%%

%%%%%%%%%%%%%%%%%%%%%%%%%%%%%%%%%%%%%%%%%%%%%%%%%%%%%%%%%%%%%%%%%%%%%%%%%%%%%%%%%%%%%%%%%%%%%%%%

%%%%%%%%%%%%%%%%%%%%%%%%%%%%%%%%%%%%%%%%%%%%%%%%%%%%%%%%%%%%%%%%%%%%%%%%%%%%%%%%%%%%%%%%%%%%%%%%

Here we show that the $L^\infty$ norm of the system \eqref{ec1d} decreases in time in the stable case ($\rho_2>\rho_1$). We will consider the set $\Omega$ equal to $\R$ or $\T$. The following theorem is the main result of the section.

\begin{thm} Let $f_0\in H^{k}(\Omega)$ with $k\geq 3$ and $\rho_2>\rho_1$. Then the unique solution to the system \eqref{ec1d} satisfies the following inequality:
$$\|f\|_{L^{\infty}}(t)\leq \|f_0\|_{L^{\infty}}.$$
\end{thm}
Proof: For $f_0\in H^{k}$ with $k\geq 3$, we prove in \cite{DY} that there exists a time $T>0$ so that the unique solution
to \eqref{ec1d} $f(x,t)$ belongs to $C^1([0,T];H^k)$.  In particular we have $f(x,t)\in C^1([0,T]\times\Omega)$, so that the Rademacher theorem gives the functions $$M(t)=\max_x f(x,t),$$ and $$m(t)=\min_x f(x,t),$$ differentiable almost every $t$. In the non periodic case, we notice that also there always exists a point $x_t\in\R$ where
$$|f(x_t,t)|=\max_x |f(x,t)|,$$ due to the fact that $f(\cdot,t)\in H^s$ with $s>1/2$, and using the
Riemann-Lebesge lemma, $f(x,t)$ tends to $0$ when $|x|\rightarrow\infty$. We suppose that this point
$x_t$ satisfies that $0<f(x_t,t)=M(t)$. A similar argument can be used for $m(t)=f(x_t,t)<0$. If we consider 
a point in which $M(t)$ is differentiable, we have
\begin{align*}
M'(t)&=\lim_{h\rightarrow 0^+}\frac{M(t+h)-M(t)}{h}\\
&=\lim_{h\rightarrow 0^+}\frac{f(x_{t+h},t+h)-f(x_t,t)}{h}\\
&=\lim_{h\rightarrow 0^+}\frac{f(x_{t+h},t+h)-f(x_t,t+h)}{h}+\frac{f(x_t,t+h)-f(x_t,t)}{h}.
\end{align*}
Since $f(x,t+h)$ takes on its maximum value at $x=x_{t+h}$, it follows:
\begin{align*}
M'(t)\!&\geq\!\lim_{h\rightarrow 0^+}\frac{f(x_t,t+h)-f(x_t,t)}{h}=f_t(x_t,t).
\end{align*}
Computing for $h>0$
\begin{align*}
M'(t)&=\lim_{h\rightarrow 0^+}\frac{M(t)-M(t-h)}{h}\\
&=\lim_{h\rightarrow 0^+}\frac{f(x_t,t)-f(x_{t-h},t-h)}{h}\\
&=\lim_{h\rightarrow 0^+}\frac{f(x_t,t-h)-f(x_{t-h},t-h)}{h}+\frac{f(x_t,t)-f(x_t,t-h)}{h}\\
&\leq \lim_{h\rightarrow 0^+} \frac{f(x_t,t)-f(x_t,t-h)}{h}\\
&\leq f_t(x_t,t),
\end{align*}
and we obtain finally
\begin{equation}
M'(t)=f_t(x_t,t).
\end{equation}
If we take the value $x=x_t$ in the equation \eqref{ec1d}, the above identity gives
\begin{equation*}
M'(t)=-\frac{\rho_2-\rho_1}{2\pi}PV\int_\R\frac{\partial_x f(x_t-\al,t)\al}{\al^2+((f(x_t,t)-f(x_t-\al,t)))^2}d\al,
\end{equation*}
using the fact that $\partial_x f(x_t,t)=0$. Integrating by parts
\begin{align*}
M'(t)&=-\frac{\rho_2-\rho_1}{2\pi}PV \int_\R\frac{\partial_{\al}(f(x_t,t)-
f(x_t-\al,t))}{\al}\D\frac{1}{1+\D\Big(\frac{f(x_t,t)-f(x_t-\al,t)}{\al}\Big)^2}d\al\\
&=I_1+I_2,
\end{align*}
where
$$I_1=-\frac{\rho_2-\rho_1}{2\pi}PV \int_\R\frac{f(x_t,t)-
f(x_t-\al,t)}{\al^2}\D\frac{1}{1+\D\Big(\frac{f(x_t,t)-f(x_t-\al,t)}{\al}\Big)^2}d\al,$$ and
$$
I_2=-\frac{\rho_2-\rho_1}{2\pi}\int_\R2\frac{\D\Big(\frac{f(x_t,t)-
f(x_t-\al,t)}{\al}\Big)^2}{\Big(1+\D\Big(\frac{f(x_t,t)-f(x_t-\al,t)}{\al}\Big)^2\Big)^2}\partial_\al
\Big(\frac{f(x_t,t)-f(x_t-\al,t)}{\al}\Big)d\al.
$$
Using the function 
$$
G(x)=-\frac{x}{1+x^2}+\arctan x,
$$
we can write $I_2$ as follows:
$$
I_2=-\frac{\rho_2-\rho_1}{2\pi}PV \int_\R \partial_\al\, G\Big(\frac{f(x_t,t)-f(x_t-\al,t)}{\al}\Big) d\al.
$$
Integrating we obtain
$$
I_2=-\frac{\rho_2-\rho_1}{2\pi}[G\Big(\lim_{\al\rightarrow+\infty}\frac{f(x_t,t)-f(x_t-\al,t)}{\al}\Big)-G\Big(\lim_{\al\rightarrow-\infty}\frac{f(x_t,t)-f(x_t-\al,t)}{\al}\Big)]=0.
$$
The $I_1$ term is equal to
$$
I_1=-\frac{\rho_2-\rho_1}{2\pi}PV\int_\R\frac{M(t)-f(x_t-\al,t)}{\al^2+(M(t)-f(x_t-\al,t))^2}d\al\leq 0,
$$
so that $M'(t)\leq 0$ for almost every $t$. In a similar way we obtain for $m(t)$ the following:
\begin{equation*}
m'(t)=-\frac{\rho_2-\rho_1}{2\pi}PV\int_\R
\frac{m(t)-f(x_t-\al,t)}{\alpha^2+(m(t)-f(x_t-\al,t))^2}d\al\geq 0,
\end{equation*}
for almost every $t$. Integrating in time we conclude the argument obtaining the maximum
principle.

%%%%%%%%%%%%%%%%%%%%%%%%%%%%%%%%%%%%%%%%%%%%%%%%%%%%%%%%%%%%%%%%%%%%%%%%%%%%%%%%%%%%%%%%%%%%%%%

Let $\Omega=\T$, with this maximum principle we can conclude the following decay of the $L^{\infty}$ norm.

\begin{prop} Let $f_0\in H^{k}(\T)$ with $k\geq 3$ and $\rho_2>\rho_1$. If 
$$\int_\T f_0(x) dx=0,$$
then the unique solution to the system \eqref{ec1d} satisfies the following inequality:
$$\|f\|_{L^{\infty}}(t)\leq \|f_0\|_{L^{\infty}}e^{-(\rho_2-\rho_1)C(\|f_0\|_{L^{\infty}})t},$$
with $C(\|f_0\|_{L^{\infty}})>0.$
\end{prop}
Proof: Suppose that 
$$\int_\T f_0(x) dx=0.$$ We can write \eqref{ec1d} as follows:
\begin{align*}
\begin{split}
\D f_t(x,t)&=\frac{\rho_2-\rho_1}{2\pi}PV\int_{\R}\partial_x \arctan\Big(\frac{f(x,t)-
f(x-\al,t)}{\al}\Big)d\al,
\end{split}
\end{align*}
and therefore we have
\begin{align*}
\begin{split}
\D \int_{\T} f_t(x,t) dx &=\frac{\rho_2-\rho_1}{2\pi}\int_\T PV\int_{\R}\partial_x \arctan\Big(\frac{f(x,t)-
f(x-\al,t)}{\al}\Big)d\al dx\\
&=\frac{\rho_2-\rho_1}{2\pi}PV \int_{\R}\int_\T \partial_x \arctan\Big(\frac{f(x,t)-
f(x-\al,t)}{\al}\Big) dx d\al\\
&=0.
\end{split}
\end{align*}
Integrating in time we obtain
\begin{equation}\label{meanzero}
 \int_\T f(x,t) dx=0, \quad \forall t\geq 0.
\end{equation}
As we show in the proof of the above theorem, we have
$$
\frac{d}{dt}
\|f\|_{L^{\infty}}(t)=-\frac{\rho_2-\rho_1}{2\pi}PV\int_\R\frac{\|f\|_{L^{\infty}}(t)-f(x_t-\al,t)}{\al^2+
(\|f\|_{L^{\infty}}(t)-f(x_t-\al,t))^2}d\al,
$$ for almost every $t$.
Using the maximum principle, for $|\al|\leq r$ we get 
$$\al^2+(\|f\|_{L^{\infty}}(t)-f(x_t-\al,t))^2\leq r^2+4\|f_0\|^2_{L^{\infty}},$$
and it follows:
\begin{align*}
\begin{split}
\frac{d}{dt}
\|f\|_{L^{\infty}}(t)&\leq-\frac{\rho_2-\rho_1}{2\pi}PV\int_{|\al|\leq r}\frac{\|f\|_{L^{\infty}}(t)-f(x_t-\al,t)}{\al^2+
(\|f\|_{L^{\infty}}(t)-f(x_t-\al,t))^2}d\al\\
&\leq -\frac{\rho_2-\rho_1}{2\pi}\frac{2r}{r^2+
4\|f_0\|^2_{L^{\infty}}}\|f\|_{L^{\infty}}(t)+\frac{\rho_2-\rho_1}{2\pi}\frac{1}{r^2+
4\|f_0\|^2_{L^{\infty}}}\int_{|\al|\leq r}f(x_t-\al) d\al.
\end{split}
\end{align*}
If we take $r=n\pi$ for $n\in\N$, using \eqref{meanzero} we obtain
$$
\frac{d}{dt}
\|f\|_{L^{\infty}}(t)\leq -\frac{\rho_2-\rho_1}{2\pi}\frac{2n\pi}{n^2\pi^2+
4\|f_0\|^2_{L^{\infty}}}\|f\|_{L^{\infty}}(t),
$$
and integrating in time we conclude the proof.
%$$
%\|f\|_{L^{\infty}}(t)\leq \|f_0\|_{L^{\infty}}e^{-(\rho_2-\rho_1)C(\|f_0\|_{L^{\infty}})t}.
%$$

For $\Omega=\R$ we obtain the following result.
\begin{prop} Let $f_0\in H^{k}(\R)$ with $k\geq 3$ and $\rho_2>\rho_1$. If $f_0(x)\leq 0$ or $f_0(x)\geq 0$,
then the unique solution to the system \eqref{ec1d} satisfies the following inequality:
$$\|f\|_{L^{\infty}}(t)\leq \frac{\|f_0\|_{L^{\infty}}}{1+(\rho_2-\rho_1)C(\|f_0\|_{L^{\infty}},\|f_0\|_{L^1})t} ,$$
with $C(\|f_0\|_{L^{\infty}},\|f_0\|_{L^1})>0.$
\end{prop}
Proof: Let consider $f_0(x)\geq 0$. The argument is similar in the other case. Our maximum principle shows that
\begin{equation*}
m'(t)=-\frac{\rho_2-\rho_1}{2\pi}PV\int_\R
\frac{m(t)-f(x_t-\al,t)}{\alpha^2+(m(t)-f(x_t-\al,t))^2}d\al\geq 0,
\end{equation*}
for almost every $t$, so that if $f_0(x)\geq 0$, then $f(x,t)\geq 0$. In a similar way as in the previous result, we can conclude that
$$
\D \int_{\R} f_t(x,t) dx=0,
$$
and therefore $$\int_{\R} f(x,t) dx=\int_{\R} f_0(x) dx.$$ Due to $f$ is nonnegative, we control the $L^1$ norm of the solution, so that $\|f\|_{L^1}(t)= \|f_0\|_{L^1}.$ We have $\|f\|_{L^{\infty}}(t)=f(x_t,t)$, and $$\frac{d}{dt}
\|f\|_{L^{\infty}}(t)=-I,$$ with
$$
I=\frac{\rho_2-\rho_1}{2\pi}PV\int_\R\frac{f(x_t,t)-f(x_t-\al,t)}{\al^2+
(f(x_t,t)-f(x_t-\al,t))^2}d\al,
$$ for almost every $t$. If we consider the interval $[-r,r]$ for $r>0$, 
$$U_1=\{\al\in[-r,r]:f(x_t,t)-f(x_t-\al,t)\geq f(x_t,t)/2 \},$$
and 
$$U_2=\{\al\in[-r,r]:f(x_t,t)-f(x_t-\al,t)< f(x_t,t)/2 \},$$
we get
$$
I\geq \frac{\rho_2-\rho_1}{2\pi}PV\int_{U_1}\frac{f(x_t,t)-f(x_t-\al,t)}{\al^2+
(f(x_t,t)-f(x_t-\al,t))^2}d\al\geq \frac{\rho_2-\rho_1}{2\pi}\frac{f(x_t,t)/2}{r^2+
4\|f_0\|^2_{L^{\infty}}}|U_1|.
$$
In order to estimate $|U_1|$, we use that $|U_1|=2r-|U_2|$, and
$$
\|f_0\|_{L^1}= \int_\R f(x_t-\al,t)d\al \geq \int_{U_2} f(x_t-\al,t) d\al \geq \frac{f(x_t,t)}{2}|U_2|,
$$
and therefore $|U_1|\geq 2(r-\|f_0\|_{L^1}/f(x_t,t)).$ This estimate gives
$$
I\geq \frac{\rho_2-\rho_1}{2\pi}\frac{f(x_t,t)/2}{r^2+
4\|f_0\|^2_{L^{\infty}}}|U_1|\geq\frac{\rho_2-\rho_1}{2\pi}\frac{r f(x_t,t)-\|f_0\|_{L^1}}{r^2+
4\|f_0\|^2_{L^{\infty}}},
$$ and this function reaches its maximum at $$r=\big(\|f_0\|_{L^1}+\sqrt{\|f_0\|^2_{L^1}+4\|f_0\|^2_{L^{\infty}}f^2(x_t,t)}\,\big)/f(x_t,t).$$ Using the maximum principle, easily we get
$$
I\geq \frac{\rho_2-\rho_1}{8\pi}\frac{\|f_0\|_{L^1}f^2(x_t,t)}{\|f_0\|^2_{L^1}+2\|f_0\|_{L^1}\|f_0\|^2_{L^\infty}+2\|f_0\|^4_{L^\infty}}
\geq (\rho_2-\rho_1)C(\|f_0\|_{L^1},\|f_0\|_{L^\infty})f^2(x_t,t).
$$
Finally, we obtain
$$
\frac{d}{dt}
\|f\|_{L^{\infty}}(t)\leq -(\rho_2-\rho_1)C(\|f_0\|_{L^1},\|f_0\|_{L^\infty})\|f\|^2_{L^{\infty}}(t),
$$
and integrating we end the proof.
%%%%%%%%%%%%%%%%%%%%%%%%%%%%%%%%%%%%%%%%%%%%%%%%%%%%%%%%%%%%%%%%%%%%%%%%%%%%%%%%%%%%%%%%%%%%%%%

%%%%%%%%%%%%%%%%%%%%%%%%%%%%%%%%%%%%%%%%%%%%%%%%%%%%%%%%%%%%%%%%%%%%%%%%%%%%%%%%%%%%%%%%%%%%%%%

%%%%%%%%%%%%%%%%%%%%%%%%%%%%%%%%%%%%%%%%%%%%%%%%%%%%%%%%%%%%%%%%%%%%%%%%%%%%%%%%%%%%%%%%%%%%%%%

\section{Three dimensional case (2-D interface)}

%%%%%%%%%%%%%%%%%%%%%%%%%%%%%%%%%%%%%%%%%%%%%%%%%%%%%%%%%%%%%%%%%%%%%%%%%%%%%%%%%%%%%%%%%%%%%%%

%%%%%%%%%%%%%%%%%%%%%%%%%%%%%%%%%%%%%%%%%%%%%%%%%%%%%%%%%%%%%%%%%%%%%%%%%%%%%%%%%%%%%%%%%%%%%%%

%%%%%%%%%%%%%%%%%%%%%%%%%%%%%%%%%%%%%%%%%%%%%%%%%%%%%%%%%%%%%%%%%%%%%%%%%%%%%%%%%%%%%%%%%%%%%%%

In this section, by using the same technique, we extend the maximum principle for the three dimensional stable case. We consider the set $\Omega$ the plane or the periodic setting.

\begin{thm}
Let $f_0\in H^k(\Omega)$ for $k\geq 4$, and $\rho_2>\rho_1$. Then the unique solution to \eqref{ec}
satisfies that $$\|f\|_{L^{\infty}}(t)\leq\|f_0\|_{L^{\infty}}.$$
\end{thm}

Proof: As we prove in \cite{DY}, there exists a time $T>0$ and
a unique solution $f(x,t)\in C^1([0,T];H^{k}(\Omega))$ solution of \eqref{ec}. In particular
$f(x,t)\in C^1([0,T]\times\Omega)$ using Sobolev inequalities. 
In the case $\Omega=\R^2$, there always exists a point $x_t\in \R^2$ where
$|f(x,t)|$ reaches its maximum due to the fact that $f(\cdot,t)\in H^s$ with $s>1$, and using the
Riemann-Lebesgue lemma $f(x,t)$ tends to 0 when $|x|\rightarrow\infty$. Suppose that this point is
for $M(t)=f(x_t,t)>0$. A similar argument can be used for $m(t)=f(x_t,t)<0$. By using the H.
Rademacher theorem, the function $M(t)$ is differentiable almost everywhere, and computing as before we obtain
\begin{equation}
M'(t)=f_t(x_t,t),
\end{equation}
for almost every $t.$ Using equation \eqref{ec}, the fact that $\grad f(x_t,t)=0$, and the last identity, we have
\begin{align*}
\begin{split}
M'(t)&=\frac{\rho_2-\rho_1}{4\pi}PV\int_{\R^2}\frac{-\grad
f(y,t)\cdot (x_t-y)}{[|x_t-y|^2+(f(x_t,t)-f(y,t))^2]^{3/2}}dy.\\
\end{split}
\end{align*}
Integrating by parts
\begin{align*}
\begin{split}
M'(t)&=\frac{\rho_2-\rho_1}{4\pi}PV\!\!\int_{\R^2}\!\!\grad_y (f(x_t,t)-
f(y,t))\cdot \frac{x_t-y}{|x_t-y|^3}\D\Big(1+\Big(\frac{f(x_t,t)-f(y,t)}{|x_t-y|}\Big)^2\Big)^{-3/2}\!\!dy\\
&=-\frac{\rho_2-\rho_1}{4\pi}PV\!\!\int_{\R^2}\!\!(f(x_t,t)-
f(y,t))\big(\di_y \frac{x_t-y}{|x_t-y|^3}\big)\D\Big(1+\Big(\frac{f(x_t,t)-f(y,t)}{|x_t-y|}\Big)^2\Big)^{-3/2}\!\!dy\\
&\quad -\frac{\rho_2-\rho_1}{4\pi}PV\!\!\int_{\R^2}\!\!\frac{f(x_t,t)-
f(y,t)}{|x_t-y|}\frac{x_t-y}{|x_t-y|^2}\cdot\grad_y \Big( 1+\Big(\frac{f(x_t,t)-f(y,t)}{|x_t-y|}\Big)^2\Big)^{-3/2}dy\\
&=J_1+J_2.
\end{split}
\end{align*}
We have
$$
J_2=-\frac{\rho_2-\rho_1}{4\pi}PV\int_{\R^2}\grad_y(\ln|x_t-y|)\cdot\grad_y
H\Big(\frac{f(x_t,t)-f(y,t)}{|x_t-y|}\Big)dy,
$$
where
$$H(x)=\frac{x^3}{(1+x^2)^{3/2}}.$$
The identity $\Delta_y(\ln|x_t-y|)/4\pi=\delta(x_t)$, and the following limit:
$$\lim_{y\rightarrow x_t}\frac{f(x_t,t)-f(y,t)}{|x_t-y|}=
\lim_{y\rightarrow x_t}\frac{f(x_t,t)-f(y,t)-\grad f(x_t,t)\cdot(x_t-y)}{|x_t-y|}=0,
$$
show that integrating by parts in $J_2$, we obtain
$$
J_2=\frac{\rho_2-\rho_1}{4\pi}PV\int_{\R^2}\Delta_y(\ln|x_t-y|)
H\Big(\frac{f(x_t,t)-f(y,t)}{|x_t-y|}\Big)dy=(\rho_2-\rho_1)H(0),
$$ and therefore $J_2=0$. The $J_1$ term is equal to
$$
J_1=-\frac{\rho_2-\rho_1}{4\pi}PV\int_{\R^2}\frac{M(t)-
f(y,t)}{[|x_t-y|^2+(M(t)-f(y,t))^2]^{3/2}}dy\leq0,
$$ so that $M'(t)\leq 0$ for almost every $t$. For $m(t)$ we have $m'(t)\geq 0$. 
Integrating in time we conclude the proof.

%%%%%%%%%%%%%%%%%%%%%%%%%%%%%%%%%%%%%%%%%%%%%%%%%%%%%%%%%%%%%%%%%%%%%%%%%%%%%%%%%%%%%%%%%%%%%%%%%%%%%%%%%%%%%%%%%%%%%%%%%%

As in the previous section, using this maximum principle we get the following decay of the $L^{\infty}$ norm.

\begin{prop} Let $f_0\in H^{k}(\T^2)$ with $k\geq 4$ and $\rho_2>\rho_1$. If 
$$\int_{\T^2} f_0(x) dx=0,$$
then the unique solution to the system \eqref{ec1d} satisfies the following inequality:
$$\|f\|_{L^{\infty}}(t)\leq \|f_0\|_{L^{\infty}}e^{-(\rho_2-\rho_1)C(\|f_0\|_{L^{\infty}})t},$$
with $C(\|f_0\|_{L^{\infty}})>0$
\end{prop}
Proof: We can write \eqref{ec} as follows:

\begin{align*}
\D f_t (x,t)&=\frac{\rho_2-\rho_1}{4\pi}PV\int_{\R^2}\frac{y}{|y|^2}\cdot \grad_x P\Big(\frac{f(x)-f(x-y)}{|y|}\Big)dy,\\
f(x,0)&=f_0(x),
\end{align*}
with $$P(x)=\frac{x}{\sqrt{1+x^2}}.$$
As in the previous section, checking the evolution of the integral of $f$ on $\T^2$,  we obtain
\begin{equation}\label{md}
\int_{\T^2} f(x,t) dx=0.
\end{equation}
The proof in the above theorem shows that 
$$
\frac{d}{dt}
\|f\|_{L^{\infty}}(t)=-\frac{\rho_2-\rho_1}{4\pi}PV\int_{\R^2}\frac{\|f\|_{L^{\infty}}(t)
-f(y,t)}{[|x_t-y|^2+(\|f\|_{L^{\infty}}(t)-f(y,t))^2]^{3/2}}dy,
$$ for almost every $t$. If we consider $x_t-y\in [-n\pi,n\pi]\times[-n\pi,n\pi]=A_n,$ with $n\in\N,$
we have 
$$|x_t-y|^2+(\|f\|_{L^{\infty}}(t)-f(x_t-\al,t))^2\leq 2(n\pi)^2+4\|f_0\|^2_{L^{\infty}}.$$
Using \eqref{md}, the above estimate gives
\begin{align*}
\begin{split}
\frac{d}{dt}
\|f\|_{L^{\infty}}(t)&\leq-\frac{\rho_2-\rho_1}{4\pi}PV\int_{(x_t-y)\in A_n}\frac{\|f\|_{L^{\infty}}(t)
-f(y,t)}{[|x_t-y|^2+(\|f\|_{L^{\infty}}(t)-f(y,t))^2]^{3/2}}dy\\
&\leq -\frac{\rho_2-\rho_1}{4\pi}\frac{(2n\pi)^2}{[2(n\pi)^2+4\|f_0\|^2_{L^{\infty}}]^{3/2}}\|f\|_{L^{\infty}}(t).
\end{split}
\end{align*}
Integrating in time we finish the proof.

\begin{prop} Let $f_0\in H^{k}(\R^2)$ with $k\geq 4$ and $\rho_2>\rho_1$. If $f_0(x)\leq 0$ or $f_0(x)\geq 0$,
then the unique solution to the system \eqref{ec1d} satisfies the following inequality:
$$\|f\|_{L^{\infty}}(t)\leq \frac{\|f_0\|_{L^{\infty}}}{(1+(\rho_2-\rho_1)C(\|f_0\|_{L^{\infty}},\|f_0\|_{L^1})t)^2} ,$$
with $C(\|f_0\|_{L^{\infty}},\|f_0\|_{L^1})>0.$
\end{prop}
Proof: Let us consider $f_0(x)\geq 0$. The same estimate is obtained for $f_0(x)\leq 0.$ We obtain as before $f(x,t)\geq 0$, and $\|f\|_{L^1}(t)=\|f_0\|_{L^1}.$ We have $\|f\|_{L^{\infty}}(t)=f(x_t,t)$, and $$\frac{d}{dt}\|f\|_{L^{\infty}}(t)=-J,$$ for almost every $t$, with
$$
J=\frac{\rho_2-\rho_1}{4\pi}PV\int_{\R^2}\frac{f(x_t,t)
-f(y,t)}{[|x_t-y|^2+(f(x_t,t)-f(y,t))^2]^{3/2}}dy.
$$  If we define the set $B_r(x_t)=\{y:|x_t-y|\leq r\}$ for $r>0$, 
$$V_1=\{y\in B_r(x_t):f(x_t,t)-f(y,t)\geq f(x_t,t)/2 \},$$
and 
$$V_2=\{y\in B_r(x_t):f(x_t,t)-f(y,t)< f(x_t,t)/2 \},$$
we get
$$
J\geq \frac{\rho_2-\rho_1}{4\pi}\frac{f(x_t,t)/2}{[r^2+
4\|f_0\|^2_{L^{\infty}}]^{3/2}}|V_1|.
$$
Using that $|V_1|=\pi r^2-|V_2|$, and
$$
\|f_0\|_{L^1}\geq \int_{V_2} f(y,t) dy \geq \frac{f(x_t,t)}{2}|V_2|,
$$
this estimate follows: $|V_1|\geq \pi r^2-2\|f_0\|_{L^1}/f(x_t,t).$ We have
$$
J\geq \frac{\rho_2-\rho_1}{8\pi}\frac{\pi r^2 f(x_t,t)-2\|f_0\|_{L^1}}{[r^2+
4\|f_0\|^2_{L^{\infty}}]^{3/2}},
$$ and taking
$$
r=\Big(\frac{2\|f_0\|_{L^1}/\pi+1}{f(x_t,t)}\Big)^{1/2},
$$
we find
$$
J\geq \frac{\rho_2\!-\!\rho_1}{8\pi}\frac{\pi (f(x_t,t))^{3/2}}{[1+2\|f_0\|_{L^1}/\pi+
4\|f_0\|^2_{L^{\infty}}f(x_t,t)]^{3/2}}\geq \frac{\rho_2\!-\!\rho_1}{8}\frac{(f(x_t,t))^{3/2}}{[1+2\|f_0\|_{L^1}/\pi+
4\|f_0\|^3_{L^{\infty}}]^{3/2}}.
$$
Finally, the following estimate is obtained:
$$
\frac{d}{dt}\|f\|_{L^{\infty}}(t)\leq -(\rho_2-\rho_1)C(\|f_0\|_{L^1},\|f_0\|_{L^\infty})\|f\|^{3/2}_{L^{\infty}}(t),
$$
and integrating we end the proof.

%%%%%%%%%%%%%%%%%%%%%%%%%%%%%%%%%%%%%%%%%%%%%%%%%%%%%%%%%%%%%%%%%%%%%%%%%%%%%%%%%%%%%%%%%%%%%%%%%%%%%%%%%%%%%%%%%%%%

%%%%%%%%%%%%%%%%%%%%%%%%%%%%%%%%%%%%%%%%%%%%%%%%%%%%%%%%%%%%%%%%%%%%%%%%%%%%%%%%%%%%%%%%%%%%%%%%%%%%%%%%%%%%%%%%%%%%

%%%%%%%%%%%%%%%%%%%%%%%%%%%%%%%%%%%%%%%%%%%%%%%%%%%%%%%%%%%%%%%%%%%%%%%%%%%%%%%%%%%%%%%%%%%%%%%%%%%%%%%%%%%%%%%%%%%%

\section{Small initial data}

%%%%%%%%%%%%%%%%%%%%%%%%%%%%%%%%%%%%%%%%%%%%%%%%%%%%%%%%%%%%%%%%%%%%%%%%%%%%%%%%%%%%%%%%%%%%%%%%%%%%%%%%%%%%%%%%%%%%

%%%%%%%%%%%%%%%%%%%%%%%%%%%%%%%%%%%%%%%%%%%%%%%%%%%%%%%%%%%%%%%%%%%%%%%%%%%%%%%%%%%%%%%%%%%%%%%%%%%%%%%%%%%%%%%%%%%%

%%%%%%%%%%%%%%%%%%%%%%%%%%%%%%%%%%%%%%%%%%%%%%%%%%%%%%%%%%%%%%%%%%%%%%%%%%%%%%%%%%%%%%%%%%%%%%%%%%%%%%%%%%%%%%%%%%%%

In the two-dimensional case, we prove in \cite{DY} that if the following quantity of the initial data is small
$$
\sum |\xi||\hat{f}(\xi)|,
$$ 
then there is global-in-time solution of the system \eqref{ec1d}. Here we show that if initially the $L^\infty$ norm of the first derivative is less than one, it continues less than one for all time.

\begin{lemma}
Let $f_0\in H^s$ with $s\geq 3,$ and $\|\partial_x f_0\|_{L^{\infty}}\leq 1.$ Then the unique solution of the system \eqref{ec1d} satisfies
$$
\|\partial_x f\|_{L^{\infty}}(t)<1.
$$ 
\end{lemma}
Proof: If we consider the following term in \eqref{ec1d}:

\begin{align*}
K=-\frac{\rho_2-\rho_1}{2\pi}PV\int_{\R}\frac{\dx
f(x-\al,t)\al}{\al^2+(f(x,t)-f(x-\al,t))^2}d\al,
\end{align*}
we can integrate by parts and get
\begin{align*}
K&=-\frac{\rho_2-\rho_1}{2\pi}PV\int_{\R}\frac{\partial_\al(f(x,t)-f(x-\al,t))}{\al}\D\frac{1}{1
+\D\Big(\frac{f(x,t)-f(x-\al,t)}{\al}\Big)^2}d\al\\
&=-\frac{\rho_2-\rho_1}{2\pi}PV\int_{\R}\frac{f(x,t)-f(x-\al,t)}{\al^2}\D\frac{1}{1
+\D\Big(\frac{f(x,t)-f(x-\al,t)}{\al}\Big)^2}d\al\\
&\quad-\frac{\rho_2-\rho_1}{2\pi}\int_\R2\frac{\D\Big(\frac{f(x,t)-
f(x-\al,t)}{\al}\Big)^2}{\Big(1+\D\Big(\frac{f(x,t)-f(x-\al,t)}{\al}\Big)^2\Big)^2}\partial_\al
\Big(\frac{f(x,t)-f(x-\al,t)}{\al}\Big)d\al\\
&=L_1+L_2.
\end{align*}
As we showed before
$$
L_2=-\frac{\rho_2-\rho_1}{2\pi}PV \int_\R \partial_\al\, G\Big(\frac{f(x,t)-f(x-\al,t)}{\al}\Big) d\al=0.
$$
so $K=L_2$. Making a change of variables we find the following equivalent system:
$$
f_t(x,t)=\frac{\rho_2-\rho_1}{2\pi}PV\int_\R\frac{\dx f(x,t)(x-\alpha)-
(f(x,t)-f(\al,t))}{(x-\alpha)^2+(f(x,t)-f(\al,t))^2}d\al.
$$
Taking one derivative in this formula, we have
\begin{equation}\label{pdf}
\dx f_t(x)= N_1(x)+N_2(x),
\end{equation} with
$$
N_1(x)=\frac{\rho_2-\rho_1}{2\pi} PV\!\!\int_\R\frac{\dx^2
f(x)(x\!-\!\al)}{(x\!-\!\alpha)^2\!+\!(f(x)\!-\!f(\al))^2}d\al,
$$
$$
N_2(x)=-\frac{\rho_2-\rho_1}{2\pi} PV\int_\R\frac{\dx f(x)-\triangle_\al f(x)}{(x-\alpha)^2}P(x,\al)d\al,\\
$$
where
$$
Q(x,\al)=2\,\frac{1+\dx f(x)\triangle_\al f(x)}{(1+(\triangle_\al f(x))^2)^2},
$$ and
$$
\triangle_{\al} f(x)=\frac{f(x)-f(\al)}{x-\al}.
$$
Now, we define in this section $$M(t)=\|\dx f\|_{L^\infty}(t),$$ then $M(t)=\max_x \dx f(x,t)=\dx f(x_t,t)$ where $x_t$
is the trajectory of the maximum. Similar conclusions are obtained for $m(t)=\min_x \dx f(x,t)$.
Using the Rademacher theorem as in the previous section, we have that $M'(t)=\dx f_t(x_t,t)$, and
$\dx^2 f(x_t,t)=0$. Therefore taking $x=x_t$ in \eqref{pdf}, we get
$$
M'(t)=N_2(x_t),
$$
due to $N_1(x_t)=0$. The inequality
$$|\triangle_\al f(x_t)|\leq M(t),$$ shows that for $M(t)<1$ the integral $N_2(x_t)\leq 0$, and therefore $M'(t)\leq 0$.
If $M(0)<1$, using the theorem of local existence, we have that for short time $M(t)<1$, and therefore
$M'(t)\leq 0$ for almost every time. This implies that $M(t)< 1$. For $m(t)$ we find $m(t)> 1.$

\subsection*{{\bf Acknowledgements}}

\smallskip

The authors were partially supported by the grant {\sc MTM2005-05980} of the MEC (Spain) and
the grant PAC-05-005-2 of the JCLM (Spain).

\begin{quote}
\begin{tabular}{ll}
Diego C\'ordoba &  Francisco Gancedo\\
{\small Instituto de Ciencias Matem\'aticas} & {\small Department of Mathematics}\\
{\small Consejo Superior de Investigaciones Cient\'ificas} & {\small University of Chicago}\\
{\small Serrano 123, 28006 Madrid, Spain} & {\small 5734 University Avenue, Chicago, IL 60637}\\
{\small Email: dcg@imaff.cfmac.csic.es} & {\small Email: fgancedo@math.uchicago.edu}
\end{tabular}
\end{quote}

%%%%%%%%%%%%%%%%%%%%%%%%%%%%%%%%%%%%%%%%%%%%%%%%%%%%%%%%%%%%%%%%%%%%%%%%%%%%%%%%%%%

\end{document}